\documentclass[preprint,12pt]{elsarticle}
\usepackage{bbding}
\usepackage{dcolumn}
\usepackage{graphicx}
\usepackage{amsmath}
\usepackage{amsfonts}
\usepackage{amssymb}
\usepackage{psfrag}
\usepackage{wrapfig}
\usepackage{subfigure}
\usepackage{makeidx}
\usepackage{bm}
\usepackage{epsf}
\usepackage{epsfig}
\usepackage{leftidx}
\usepackage{float}
\usepackage{extarrows}
\usepackage{color}
\newdefinition{rmk}{Remark}
\begin{document}

\title{An Alternative Analysis of Discontinuous Galerkin Method for Hyperbolic Conservation Law}

\author[HKUST,HKUST2]{Kun Xu\corref{cor1}}
\ead{makxu@ust.hk}

\author[HKUST]{Chang Liu}
\ead{cliuaa@ust.hk}

\author[HKUST]{Xiaodong Ren}
\ead{maxdren@ust.hk}


\address[HKUST]{Department of Mathematics, Hong Kong University of Science and Technology, Clear Water Bay, Kowloon, Hong Kong}

\address[HKUST2]{Department of Mechanical and Aerospace Engineering, Hong Kong University of Science and Technology, Clear Water Bay, Kowloon, Hong Kong}

\begin{abstract}
The development and application of the Discontinuous Galerkin (DG) method have attracted great attention in computational fluid dynamics (CFD) community in the past decades.
The underlying reason for such an intensive investigation is due to favorable properties of the DG method, such as higher order, compactness, and easy parallelization.
However, for the compressible flow simulations, the DG method is  also associated with unfavorable properties, such as the frequent instabilities in non-smooth regions.
Due to the finite element formulation, except the update of the cell averaged conservative flow variable,
it becomes rather difficulty to fully understand the physical mechanism in the DG method for the updates of other degrees of freedom
inside each element.
In this short note, based on the linear advection equation we are going to analyze the DG method in an alternative way.
The results obtained may be useful for the CFD community for its further development of the DG methods.

\end{abstract}
\maketitle

\section{Introduction}

The development of high-order CFD methods is important for delivering accurate  solutions with lower computational cost.
The discontinuous Galerkin (DG)  finite element method is one of the most well-known schemes which are under intensive investigation continuously in the past decades.
The DG  method was originally developed by
Reed and Hill \cite{reed} to solve the neutron transport
equation. The development of the DG method for
hyperbolic conservation laws was pioneered by
Cockburn, Shu and their collaborators in a series of
papers on the Runge-Kutta DG (RKDG) method \cite{Cockburn1,Cockburn2,Cockburn3,Cockburn4,Cockburn5}.
The explosive
interest in the DG method comes from its attractive
features. The DG method combines two favorable
features of the finite volume and finite element methods, i.e., the
gas evolution model of the Riemann solvers, and high accurate
solution expansion through polynomials within each element.
Due to the compactness of its discretization,
the degrees of freedom associated with any element are coupled
only with those of the neighboring elements. This feature makes the DG method be suitable
 for parallel computers. This compactness also
results in highly sparse matrices in a linearized
implicit time integration scheme.

The DG method is based on a weak formulation
of the governing equation, where the solution within each element is expanded in terms of polynomials.
In the DG framework,
besides the update of conservative flow variables, which is crucial for the conservation laws, the gradients (coefficients in front of expansion polynomials)
or the equivalent degrees of freedom (subcell solutions at special nodal points), are updated as well.
The equivalence between the DG and many other high-order methods has
been established in recent investigations \cite{huynh1, huynh2, wang1, wang2}. The conservation for the conservative flow variables as
a physical law is fully preserved in the DG method, which is the same as the finite volume scheme.
However,  the physical laws underlying the updates of other degrees of freedom, are not very clear.

All higher order methods need more degrees of freedom to present and
evolve the solution. The DG has done it in a compact way,
which is different from the finite volume ENO and
WENO approaches \cite{harten,liu,jiang}, where the additional degrees of freedom
are obtained through the use of large stencils and the information from other
cells. Then, what kind of mechanism have the DG methods been used to keep
and evolve more degrees of freedom compactly?
The main concern of the DG method is its remarkable instabilities, such as relatively easy breakdown in comparison with finite volume methods.
The instabilities observed in DG methods in the non-smooth region may have fundamental physical reasons instead of purely numerical discretization.
In this paper, based on the linear advection equation we are going to analyze the underlying physical mechanism for the DG method
in an alternative way, i.e., a different way from the commonly used analysis.
Hopefully, the mathematical and physical understanding  will be helpful for the further development of the DG methods.

\section{Mathematical Analysis}

Consider linear advection equation
\begin{equation}\label{linear-advection}
\begin{cases}
  &u_t+u_x=0  \quad \text{in}\ \text{I}=(0,1),\\
  &u(x,0)=u_0(x),
\end{cases}
\end{equation}
with periodic boundary condition.

Let $0=x_{1/2}<x_{3/2}<...<x_{N+1/2}=1,$ and denote $\text{I}_j=[x_{j-1/2},x_{j+1/2}],$ $\Delta x_j=x_{j+1/2}-x_{j-1/2},$ and $x_j=(x_{j+1/2}+x_{j-1/2})/2$.
Define $$V_h=\left\{v\in L^2(\text{I}):\ v|_{\text{I}_j}\in P^k(\text{I}_j)\right\},$$ and DG scheme is to find $u_h\in V_h$ s.t.
\begin{equation}\label{weak}
  \int_{\text{I}_j} \partial_t u_h v dx + u^{j+1/2} v^{-}_{j+1/2}-u^{j-1/2} v^{+}_{j-1/2}-\int_{\text{I}_j}u_hv_x dx =0,
\end{equation}
for all $v\in V_h$, and $j=1,...,N$.

\subsection{\bf{Second order scheme}}
Consider $k=1$ and choose the local basis in $\text{I}_j$ as
$$\left\{1,\frac{x-x_j}{x_{j+1/2}-x_{j-1/2}}\right\}.$$
Denote $\xi=\frac{x-x_j}{x_{j+1/2}-x_{j-1/2}}$, and function $u$ can be written as
\begin{equation}\label{u}
  u_h|_{\text{I}_j}=a^j_0(t)\cdot 1+a^j_1(t)\cdot \xi,
\end{equation}
in which $a_0(t)$ is the cell average of $u_h$.
Comparing to the Taylor series of $u_h|_{I_j}$ around $x_j$, the coefficients $a^j$ can be expressed by the derivatives of $u$ as
\begin{equation}\label{k1-a}
\begin{aligned}
  &a^j_0=u^j+\frac{1}{24}u^j_{xx}\Delta x_j^2,\\
  &a^j_1=u^j_x \Delta x_j +\frac{1}{40} u^j_{xxx}\Delta x_j^3,
\end{aligned}
\end{equation}
where $u^j$, $u^j_x$, $u^j_{xx}$, $u^j_{xxx}$ are the values of function $u$ and its derivatives at $x_j$.

Let $v=1$ in Eq.\eqref{weak}, we have
\begin{equation}\label{k1-0-1}
  \frac{d a^j_0}{dt}+\frac{1}{\Delta x_j}(u^{j+1/2}-u^{j-1/2})=0.
\end{equation}
When the flux at cell interface is chosen as the upwind flux
$$u^{j+1/2}=u^{j+1/2,-}_h, u^{j-1/2}=u^{j-1/2,-}_h,$$
 Eq.\eqref{k1-0-1} becomes the upwind scheme
\begin{equation}\label{k1-0-2}
  \frac{d a^j_0}{dt}+\frac{1}{\Delta x_j}\left\{\left(a^j_0+\frac12 a^j_1\right)-\left(a^{j-1}_0+\frac12 a_1^{j-1}\right)\right\}=0,
\end{equation}
which is basically the same as the second order finite volume scheme.

Let $v=\xi$ in Eq.\eqref{weak}, we have
\begin{equation}\label{k1-1-1}
  \frac{da^j_1}{dt}+\frac{6}{\Delta x_j}\left(u^{j+1/2}-2a^j_0+u^{j-1/2}\right)=0,
\end{equation}
which corresponds to
\begin{equation} \label{eq:32}
  u_{xt}+ u_{xx}=O(\Delta x_j^2),
\end{equation}
by expanding terms in Eq.\eqref{k1-1-1} around $x_j$.
When using upwind flux, Eq.\eqref{k1-1-1} becomes
\begin{equation}\label{k1-1-2}
  \frac{da^j_1}{dt}+\frac{6}{\Delta x_j}\left(a^j_0+\frac12 a^j_1+a^{j-1}_0+\frac12 a^{j-1}_1-2a^j_0\right)=0,
\end{equation}
which corresponds to
\begin{equation}\label{k1-1-2-t}
  u_{xt}+\frac{2}{5}u_{xxx}\Delta x_j=O(\Delta x_j^2).
\end{equation}
In such a case, the second order scheme degenerates to the first order upwind scheme. This may be the reason that the second order finite volume upwind scheme
seems have better accuracy than the second-order upwind DG method.
Eq.\eqref{k1-1-2} will be consistent to Eq.\eqref{eq:32} if a higher-order expansion for the solution is used, see next subsection.

The update of second order DG scheme with upwind flux can be written as
\begin{equation}\label{DG-2}
  \Delta x_j \frac{d}{dt}\left(\begin{array}{c}a_0^j \\a_1^j \\\end{array}\right)
  +\left(\begin{array}{cc}1 & 1/2 \\-6 & 3 \\\end{array}\right)
  \left(\begin{array}{c}a_0^j \\a_1^j \\\end{array}\right)
  -\left(\begin{array}{cc}1 & 1/2 \\-6 & -3 \\\end{array}\right)
  \left(\begin{array}{c}a_0^{j-1} \\a_1^{j-1} \\\end{array}\right)=0 .
\end{equation}

\subsection{\bf{Third order scheme}}

For $k=2$, we choose the local basis as $$\left\{1,2\sqrt{3}\xi,6\sqrt{5}\xi^2-\sqrt{5}/2\right\}$$ with $\xi=\frac{x-x_j}{x_{j+1/2}-x_{j-1/2}}$,
 and function $u$ can be written as
\begin{equation}
\begin{aligned}
  u_h|_{\text{I}_j} & =a^j_0(t)\cdot1 + a^j_1(t) \cdot(2\sqrt{3}\xi)  + a_2^j(t) \cdot(6\sqrt{5}\xi^2-\sqrt{5}/2).
\end{aligned}
\end{equation}
Comparing to the Taylor series of $u_h|_{I_j}$ around $x_j$, the coefficients $a^j$ can be expressed by the derivatives of $u$ as
\begin{equation}\label{k2-a}
\begin{aligned}
  &a^j_0=u^j+\frac{1}{24}u^j_{xx}\Delta x_j^2,\\
  &a^j_1=\frac{1}{2\sqrt{3}}u^j_x \Delta x_j +\frac{1}{80\sqrt{3}} u^j_{xxx}\Delta x_j^3,\\
  &a^j_2=\frac{1}{12\sqrt{5}}u^j_{xx} \Delta x_j^2.
\end{aligned}
\end{equation}
Let $v=1$ in Eq.\eqref{weak}, we have
\begin{equation}\label{k2-1-1}
  \frac{d a^j_0}{dt}+\frac{1}{\Delta x_j}(u^{j+1/2}-u^{j-1/2})=0.
\end{equation}
If the flux at cell interface is chosen as the upwind flux
$$u^{j+1/2}=u^{j+1/2,-}_h, u^{j-1/2}=u^{j-1/2,-}_h,$$
Eq.\eqref{k2-1-1} becomes the upwind scheme
\begin{equation}\label{k2-1-2}
  \frac{d a^j_0}{dt}+\frac{1}{\Delta x_j}\left\{\left(a^j_0+\sqrt{3}a^j_1+\sqrt{5}a^j_2\right)-\left( a^{j-1}_0+\sqrt{3}a^{j-1}_1+\sqrt{5}a^{j-1}_2\right)\right\}=0.
\end{equation}

Let $v=2\sqrt{3}\xi$ in Eq.\eqref{weak}, we have
\begin{equation}\label{k2-2-1}
  \frac{d a^j_1}{dt}+
  \frac{\sqrt{3}}{\Delta x_j}\left(u^{j+1/2}-2a^j_0+u^{j-1/2}\right)=0,
\end{equation}
which can be rewritten as a finite volume scheme plus a correction term $C\sim O(\Delta x_j^2)$
 \begin{equation}\label{k2-2-2}
   \frac{dU^j_x}{dt}+\frac{1}{\Delta x_j/2}\left(F_{i+1/4}-F_{i-1/4}\right)=C,
 \end{equation}
 where the fluxes at $x_{i+1/4}$ and $x_{i-1/4}$ are
 \begin{equation}
 \begin{aligned}
   &F_{i+1/4}=\frac{U^{j+1/2}-U^j}{\Delta x_j/2},\\
   &F_{i-1/4}=\frac{U^j-U^{j-1/2}}{\Delta x_j/2},
 \end{aligned}
 \end{equation}
 and the correction term is
 \begin{equation}\label{correction}
  \begin{aligned}
   C=&\frac{2(U^{j+1/2}+U^{j-1/2}-2U^j)}{\Delta x_j^2}-\frac12U_{xx}^j\\
   &\sim O(\Delta x_j^2).
  \end{aligned}
 \end{equation}
If we choose the cell interface flux as the upwind flux, Eq.\eqref{k2-2-1} becomes
\begin{equation}\label{k2-2-3}
  \begin{aligned}
  \frac{da^j_1}{dt}+\frac{\sqrt{3}}{\Delta x_j}&\left((a^j_0+\sqrt{3}a^j_1+\sqrt{5}a^j_2)\right.\\
  &\left.+(a^{j-1}_0+\sqrt{3}a^{j-1}_1+\sqrt{5}a^{j-1}_2)-2a^j_0\right)=0.
  \end{aligned}
\end{equation}
Both Eq.\eqref{k2-2-1} and Eq.\eqref{k2-2-3} correspond to
\begin{equation}\label{k2-2-3-t}
  u_{xt}+u_{xx}=O(\Delta x_j^2),
\end{equation}
which is the basic evolution law for the expansion coefficient $a_1^j$ or $U_x^j$.

Let $v=6\sqrt{5}\xi^2-\sqrt{5}/2$ in Eq.\eqref{weak}, we have
\begin{equation}\label{k2-3-1}
  \frac{da^j_2}{dt}+\frac{\sqrt{5}}{\Delta x_j}\left(u_h^{j+1/2}-u_h^{j-1/2}-\sqrt{12}a^j_1\right)=0.
\end{equation}
If we use upwind flux, Eq.\eqref{k2-3-1} becomes
\begin{equation}\label{k2-3-2}
\begin{aligned}
  \frac{da^j_2}{dt}+\frac{\sqrt{5}}{\Delta x_j}&\left((a^j_0+\sqrt{3}a^j_1+\sqrt{5}a^j_2)\right.\\
  &\left.-(a^{j-1}_0+\sqrt{3}a^{j-1}_1+\sqrt{5}a^{j-1}_2)-\sqrt{12}a^j_1\right)=0.
\end{aligned}
\end{equation}
Both Eq.\eqref{k2-3-1} and Eq.\eqref{k2-3-2} correspond to
\begin{equation}\label{k2-3-2-t}
  u_{xxt}+ u_{xxx}=O(\Delta x_j^2),
\end{equation}
which is the basic evolution law for the expansion coefficient $a_2^j$ or $U_{xx}^j$.

The third order DG scheme with upwind flux can be written as
\begin{equation}\label{DG-3}
\begin{aligned}
  \Delta x_j\frac{d}{dt}\left(\begin{array}{c}a_0^j \\a^j_1 \\a_2^j \\\end{array}\right) &
  +\left(\begin{array}{ccc}1&\sqrt{3}&\sqrt{5} \\-\sqrt{3} & 3 & \sqrt{15}\\ \sqrt{5}&-\sqrt{15}&5 \\\end{array}
   \right)
   \left(\begin{array}{c}a_0^j \\a^j_1 \\a_2^j \\\end{array}\right) \\
&  -\left(\begin{array}{ccc}1&\sqrt{3}&\sqrt{5} \\-\sqrt{3} & -3 & -\sqrt{15}\\\sqrt{5}&\sqrt{15}&5 \\\end{array}
   \right)
   \left(\begin{array}{c}a_0^{j-1} \\a_1^{j-1} \\a_2^{j-1} \\\end{array}\right)=0 .
   \end{aligned}
\end{equation}

Based on the above analysis, we can realize that for the compact DG method,
the updates of the expansion coefficients are basically coming from the additional spatial derivative of the conservation law  (\ref{linear-advection}),
such as Eq.\eqref{k2-2-3-t} and \eqref{k2-3-2-t}, and this process can be continued with even higher order expansions of the solution inside each element.

\section{Physical Understanding of DG Method}

The analysis in the last section figures out the underlying mechanism for the update of additional degrees of freedom in the
DG method. Based on the above analysis, we may get a better understanding of the DG method, and explain numerical observations easily.
The first concern is about the robustness of the scheme. As we know, the conservation law cannot be mathematically manipulated due to the
emerging of discontinuous solutions. So, in the discontinuous region, the updates of high order parameters in the DG method are basically unreliable.
They are based on faked evolution equations. The above analysis fully supports the correctness of the DG method for the smooth flow.
Although the DG is associated with the name "discontinuous", it is intrinsically based  on the smoothness assumption in order to validate its evolution mechanism.
This is why the limiters and filtering techniques associated with DG method play such an important role, and many efforts have been continuously devoted to
the development of this kind of treatments.
For the finite volume method, the higher-order is coming from reconstruction and the
information from a large stencil is used. There is no dynamics for the additional degrees of freedom in a finite volume approach,
 and the reconstruction takes place at the beginning of each time step and the limiter is used before the evolution of
the solution, such as a pre-processing of the initial condition.
However, for the DG method the high-order is achieved dynamically through the evolution of higher-order terms, and these terms
are governed by additional equations which are manipulated  from the original one. The limiter used in DG is a post-processing of the solution.
Since there is no a clear boundary between the well-resolved (smooth) and under-resolved (discontinuous) solutions, it should be
extremely cautious in the use of DG in comparison with the finite volume one.
It is not surprising for the observation that the DG solution is much more sensitive to limiter than the finite volume scheme.
For both finite volume and DG methods, the accuracy depends closely on the limiters on the higher-order degrees of freedom. For finite volume method,
much more information is provided through large stencils, which are very reliable through the WENO formulation. However, based on the same limiting idea,
due to compactness the DG uses a much narrow stencil of neighboring cells to limit the higher-order degrees of freedom, which may not have enough information to do it.
If the stencils for the DG method are much enlarged, the differences between finite volume and DG will diminish. So, for non-smooth flow it is contradictory in the DG formulation
for the compactness and accuracy. This is why great efforts have been continuously devoted to the design of new limiters.

Since the additional degrees of freedom in DG are updated following additional evolution equations,
the DG method will encounter difficulties if the original governing equation is much more complicated than the hyperbolic one, such as the Navier-Stokes equations or Burnett equations,
because the additional spatial derivatives will make these equations more un-realizable physically.
Therefore, it is necessary to change the complicated system into  first order equations with a much enlarged number of first order equations \cite{bassi,cockburn}.
The first order PDEs can be much easily manipulated, such as the linear advection equation.
Theoretically, it is perfectly valid to
take additional derivatives of the original PDE in the smooth flow region. That is why the DG method works perfectly in flow problem with smooth solution or other problems without discontinuities.
However, for high speed flow problems, especially with complicated geometries and wave interactions,
it will become hard to distinguish the smooth and discontinuous regions.

In order to construct a higher-order compact scheme, an accurate capturing
of both spatial and temporal higher-order evolution solution from the original governing equation may become necessary.
In this sense, for the gas dynamic equations, the first order Riemann solution has to be abandoned for a genuinely higher-order scheme.
A generalized Riemann
solver based on the initial piecewise discontinuous linear flow distribution
has been developed \cite{ben}, which shows better numerical performance
than the 1st-order Riemann-solution based methods.
 Even though it is difficult, it is fully legitimate to
continue this track and to construct schemes based on the time evolution
solution from  more general piecewise discontinuous polynomials as an initial condition. The ADER method
\cite{toro} is based on such a understanding, even though the linearized evolution
equations for the higher-order derivatives are adopted, which have the similar underlying mechanism as the DG one.

The development of higher-order scheme should be based on the higher-order gas evolution model by capturing the evolution solution of the basic conservation law
from a more general initial condition,
such as an initially discontinuous piecewise polynomials on both sides of a cell interface.
Based on a higher-order evolution model, besides updating the
conservative flow variables through the cell interface fluxes, the evolution solution also provides
the time-dependent solutions at the subcell locations or at the cell interface at the next time level.
Therefore,  more information, such as cell averaged and subcell or cell interface values, are available to
be used for the reconstruction of higher-order initial condition at the next time level compactly.
In this sense, due to higher-order evolution model, the information or stencils needed for the higher-order WENO reconstruction can be obtained locally and compactly.
For example, for the third-order compact gas-kinetic scheme \cite{pan,xu},
besides the update of the cell averaged value,
the values at the cell interface are updated as well for the initial reconstruction.
The development of such a higher-order compact scheme depends on the higher evolution
model rather than the first order Riemann solution.
The first order Riemann solution doesn't provide the
necessary dynamics for a higher-order compact scheme.
That is why the DG method must use the weak formulation, but this formulation has a weak physical basis for the
update of additional degrees of freedom.
The gas-kinetic scheme has been developed in the DG framework as well \cite{luo,ren},
which seems more robust than the standard DG method, especially for the viscous and heat conducting flows.
One of the main reasons for this is that even with a discontinuous initial condition, the evolution mechanism, such as the initial particle free transport, immediately generates a continuous solution of the gas distribution function.
At the same time, the viscosity coefficient in the gas-kinetic DG scheme (GKS-DG) is enhanced through the increment of the particle collision time at the cell interface
according to the pressure jump. As a result, the thickness of the "discontinuous layer" gets enlarged to the scale of mesh size.
Therefore, the GKS-DG method is intrinsically an artificial dissipation DG approach, and the solutions are supposed to be a continuous one.
For a smooth solution, it is legitimate to update the additional degrees of freedom through the DG framework.
But, in comparison with finite volume compact GKS, the robustness of the compact GKS-DG is inferior.

\section{Conclusion}

In this paper, based on the linear advection equation, we have analyzed the physical mechanism of the DG method for the updates of its degrees of freedom.
With the expansion of the solution inside each element, same as the finite volume scheme the DG updates the cell averaged conservative flow variable conservatively.
Due to its compactness, the updates of additional degrees of freedom in DG are based on the new governing equations, which
are obtained implicitly through the spatial derivatives of the original physical law.
Mathematically, the DG can be  considered as a higher-order Lax-Wendroff method in space, instead of the original Lax-Wendroff technique in time.

In order to construct a reliable compact higher-order method, it seems necessary to use a higher order evolution solution of the original conservation law, such as the
time dependent solution evolved from a piecewise discontinuous $P2$ distribution around a cell interface directly, the so-called "generalized" Riemann problem.
Based on such an evolution solution, besides the flux evaluation the time-dependent flow distribution at some locations inside each element or at the cell interface
can be obtained explicitly as well, which can be used directly for a higher-order local
compact initial reconstruction at the next time level. In other words, the higher-order evolution model can effectively move the large stencils needed in the finite volume WENO reconstruction to
the local neighboring cells.
The difficulty for properly limiting higher-order degrees of freedom in the DG method is due to its first order
Riemann solution, which doesn't provide enough information through the averaged conservative flow variables at the neighboring cells.
With higher-order evolution model, the weak formulation used in the DG method for the updates of additional degrees of freedom, such as the underlying governing equations presented in this
paper, can be ignored.
In other words, there is no need to adopt new governing equations.
In conclusion, the real obstacle for developing reliable compact higher-order schemes is the use of the first order Riemann solution, where its dynamics is insufficient for the
construction of higher-order compact schemes.

\subsection*{{\bf Acknowledgements}}

The current work was supported by Hong Kong Research Grant Council. We thank Professors Y. Xu, J.Q. Li, H. Xiao, and Y.X. Ren for helpful discussion, even though they may not
fully support the arguments in this paper.

\end{document}